\documentclass [11pt,twoside,a4paper]{article}
\usepackage{amsfonts}
\usepackage{amsthm}
\usepackage{amsmath}
\usepackage{amstext}
\usepackage{amssymb}
\usepackage{mathrsfs}
\usepackage{amscd}
\usepackage{subfigure}
\usepackage{mathabx}
\usepackage[all]{xy}
\usepackage{xypic}
\usepackage{epsf} 
\usepackage{graphicx} 
\usepackage{fancybox} 
\usepackage{color} 
\usepackage{fancyhdr}
\usepackage[hang,footnotesize]{caption2} 

\def\mathcal{\mathscr}
\newfont{\aaa}{cmb10 at 19pt}
\newfont{\bbb}{cmb10 at 11pt}

\newtheorem{thm}{Theorem}

\newcommand{\bm}[1]{\mbox{\boldmath $#1$}}

\def\v1{\vspace{1mm}}

\def\nnd{\noindent}  \def\lb{\label}
\def\qd{\quad}   \def\qqd{\qquad}
\def\lz{\lambda}  \def\dz{\delta}

\def\bg{\begin}
\def\be{\bg{equation}}     \def\de{\end{equation}}
\def\bgar{\bg{eqnarray}}   \def\dear{\end{eqnarray}}
\def\ben{\bg{enumerate}}    \def\den{\end{enumerate}}
\def\thm{\nnd\bg{thm1}}     \def\dethm{\end{thm1}}

\def\algo{\nnd\bg{algo1}}   \def\dealgo{\end{algo1}}

\def\xmp{\nnd\bg{xmp1}}     \def\dexmp{\end{xmp1}}

\def\lmm{\nnd\bg{lmm1}}     \def\delmm{\end{lmm1}}

\def\rmk{\nnd\bg{rmk1}}     \def\dermk{\end{rmk1}}

\def\prp{\nnd\bg{prp1}}     \def\deprp{\end{prp1}}
\def\prf{\medskip \noindent {\bf Proof}. }    \def\deprf{\qed\medskip}

\newcommand{\rf}[2]{[\ref{#1}; #2]}

\font\cms=cmss9 scaled \magstep1

\makeatletter
\def\@evenhead{
\vbox{\hbox to \textwidth {}{\hspace{0mm}{\footnotesize
\thepage}}{\hspace{8cm} {\footnotesize {Yufeng GAO, Yanxun
CHANG}}} \protect\vspace{1truemm}\relax \hrule depth0pt
height0.15truemm width\textwidth}}
\def\@evenfoot{}
\def\@oddhead{\vbox{\hbox to \textwidth
{{\hspace{0cm}{\footnotesize Symmetric $\lambda$-configurations
with small $\lambda$}\hfill{\footnotesize
\thepage}}\hspace{0mm}}{} \protect\vspace{1truemm}\relax\hrule
depth0pt height0.15truemm width\textwidth}}
\def\@oddfoot{}
\makeatother

\begin{document}

\thispagestyle{empty} \thispagestyle{fancy} {
\fancyhead[lO,RE]{\footnotesize Front. Math. China \\
https:/\!/doi.org/10.1007/s11464-018-0717-9\\[3mm]
}
\fancyhead[RO,LE]{\scriptsize \bf 
} \fancyfoot[CE,CO]{}}
\renewcommand{\headrulewidth}{0pt}


\setcounter{page}{1}
\qquad\\[8mm]

\noindent{\aaa{Approximation Theorem for Principal Eigenvalue of Discrete {\LARGE $\bm p$}-Laplacian}}\\[1mm]

\noindent{\bbb Yueshuang LI}\\[-1mm]

\noindent\footnotesize{School of Mathematical Sciences, Beijing Normal University, Beijing 100875, China}\\[6mm]

\vskip-2mm \noindent{\footnotesize$\copyright$ Higher Education
Press and Springer-Verlag GmbH Germany 2018, part of Springer Nature 2018} \vskip 4mm

\normalsize\noindent{\bbb Abstract}\quad For the principal eigenvalue of discrete weighted $p$-Laplacian on the set of nonnegative integers, the convergence of an approximation procedure and the inverse iteration is proved. Meanwhile, in the proof of the convergence,
the monotonicity of an approximation sequence is also checked.
To illustrate these results, some examples are presented.\vspace{0.3cm}

\footnotetext{Received June 8, 2018; accepted July 17,
2018\\
\hspace*{5.8mm} E-mail:
liyueshuang@mail.bnu.edu.cn}

\noindent{\bbb Keywords}\quad Principal eigenvalue, weighted~$p$-Laplacian, approximation theorem, inverse iteration\\
{\bbb MSC}\quad 34L15, 34G20, 39A12, 60J27, 65F15\\[0.4cm]

\noindent{\bbb 1\quad Introduction}\lb{s-01}\\[0.1cm]

In order to compute the principal eigenvalue of discrete weighted $p$-Laplacian, a series of results were presented in \cite{CWZ}. In fact, the results are more or less complete, including the variational formulas in different formulations, explicit lower and upper estimates, a criterion for positivity and an approximating procedure for the principal eigenvalue.
The approximating procedure is very much practical in computation. It is regretted that the convergence is still missed in that setup. This is the aim of this paper.

We need some notation to state our main results. Here, the notation and boundary conditions are the same as those used in \rf{CWZ}{\S  2}, unless otherwise stated. Let $E:=\{k\in \mathbb{Z}_+:~0\leq k < N+1\}$ with $N\leq \infty$ and $\{\nu_k: k\in E\}$ is a positive sequence with boundary condition $\nu_{-1}=0$. First, we define an operator $I\!I$ as follows
\begin{align*}
I\!I_{i}(f)=\frac{1}{f_{i}^{p-1}}\left[\sum_{j=i}^{N}\hat{\nu}_{j}
\left(\sum_{k=0}^{j}\mu_{k}f_{k}^{p-1}\right)^{p^{*}-1}\right]^{p-1},\qqd i\in E,
\end{align*}
where $f$ is a positive function on $E$, $p^{*}$ is the conjugate number of $p$ (i.e.,$1/p+1/p^{*}=1$), $\{\mu_k: k\in E\}$ is a positive sequence and $\hat{\nu}_{j}=\nu_{j}^{1-p^{*}}$.
In fact, the operator $I\!I$ in linear case was first introduced in \cite{CMF} by Chen in 2001.
Next, let $\partial_{k}(f)=f_{k+1}-f_{k},$
then the weighted $p$-Laplacian $\Omega_{p}$ ($p>1$) can be rewritten as
\begin{equation} \Omega_{p}f(k)=\nu_{k}\mid\partial_{k}(f)\mid^{p-1}\text{sgn}(\partial_{k}(f))-
\nu_{k-1}\mid\partial_{k-1}(f)\mid^{p-1}\text{sgn}(\partial_{k-1}(f)).   \lb{01}\end{equation}
Let $\lz_{p}$ denote the principal eigenvalue of the weighted $p$-Laplacian $\Omega_p$ and set
$$\sigma_{p}=\sup_{n \in E }\mu[0,n]\hat{\nu}[n,N]^{p-1},$$
where we write $\mu[m,n]=\sum_{j=m}^{n}\mu_{j}$ for a measure $\mu.$
In what follows, let $f^{(n)}_k$ denote the $k$th component of the vector $f^{(n)}$
and let $$D_{p}(f)=\sum_{k\in E}\nu_k|f_k-f_{k+1}|^p.$$
When $N\leq\infty$, the approximation procedure for $\lz_{p}$ is presented in the following theorem.

\thm\lb{t-01}\qd{\cms Assume that $N\leq\infty$  and $\sigma_{p}<\infty$.
\bg{itemize} \setlength{\itemsep}{-0.6ex}
\item[(1)] Define
$$f^{(1)}=\hat{\nu}[\cdot,N]^{1/p^*},\qd f^{(n)}=f^{(n-1)}\left(I\!I\left(f^{(n-1)}\right)\right)^{p^{*}-1},\qd n\geq 2,$$
and $$\dz_n=\sup\limits_{i\in E}I\!I_{i}\left(f^{(n)}\right).$$
Then $\dz_n$ is decreasing in $n$~(denote its limit by $\dz_{\infty}$) and
$$\lz_p\geq \dz_{\infty}^{-1}\geq\cdots\geq\dz_1^{-1}>0.$$
\item[(2)] For fixed $\ell,m \in E$,~$\ell<m,$ define
$$f^{(1,\ell,m)}=\hat{\nu}[\cdot\vee\ell,m]\mathbb{1}_{\leq m},$$
$$f^{(n,\ell,m)}=f^{(n-1,\ell,m)}\left(I\!I\left(f^{(n-1,\ell,m)}\right)\right)^{p^*-1}
\mathbb{1}_{\leq m},\qd n\geq 2,$$
where $\mathbb{1}_{\leq m}$ is the indicator of the set $\{0,1,\cdots,m\}$ and then define
$$\dz_n^{'}=\sup_{\ell,m:\ell<m}\min_{i\leq m}I\!I_i\left(f^{(n,\ell,m)}\right).$$
Then $\dz_n^{'}$ is increasing in $n$~(denote its limit by $\dz_{\infty}^{'}$)~and
$$\sigma_{p}^{-1}\geq\dz_1^{'-1}\geq\cdots\geq\dz_{\infty}^{'-1}\geq\lz_p.$$
Next, define
$$\bar{\dz}_n=\sup_{\ell,m:\ell<m}\frac{\mu\left(f^{(n,\ell,m)^p}\right)}
{D_{p}\left(f^{(n,\ell,m)}\right)},\qqd n\geq 1.$$
Then $\bar{\dz}_n^{-1}\geq\lz_p$ and $\bar{\dz}_{n+1}\geq\dz_{n}^{'}$ for $n\geq 1.$
\end{itemize}}\dethm

In fact, as an extension of \cite{Chen10}, Theorem \ref{t-01} was first presented in \rf{CWZ}{Theorem  2.4}, where the monotonicity of $\{\delta_{n}\}$ and $\{\delta_{n}'\}$ was proved, besides, the comparison property of $\{\bar{\dz}_n\}$ and $\{\dz_{n}^{'}\}$ was also proved there. Here, we remark that three questions are open in Theorem \ref{t-01}: 
\bg{itemize} \setlength{\itemsep}{-0.6ex}
\item[(i)] $\dz_{\infty}^{-1}\stackrel{?}{=}\lz_p$; 
\item[(ii)]$\dz_{\infty}^{'-1}\stackrel{?}{=}\lz_p$; 
\item[(iii)] does $\{\bar{\dz}_n\}$ have monotonicity in $n?$ 
\end{itemize}
The main aim of this paper is to answer these three questions in the case of finite $N$ and the answers are stated in the following theorem.

\thm\lb{t-02}\qd{\cms Assume that $N<\infty.$ Let $f^{(1)}>0$ be given.
\bg{itemize} \setlength{\itemsep}{-0.6ex}
\item[(1)] Define successively functions $f^{(n)}$ on $E:$
$$f^{(n)}=f^{(n-1)}\left(I\!I\left(f^{(n-1)}\right)\right)^{p^{*}-1},\qqd n\geq 2.$$
\item[(2)] For $n\geq 1$, define three sequences as follows
 $$\delta_{n}=\sup_{i \in E} I\!I_{i}\left(f^{(n)}\right),\qqd \delta_{n}'=\inf_{i\in E}I\!I_{i}\left(f^{(n)}\right),\qqd \overline{\delta}_{n}=\frac{\mu\left({f^{(n)^p}}\right)}{ D_{p}\left({f^{(n)}}\right)}.$$
Then $\{\delta_{n}\}$ is decreasing in $n$, both $\{\delta_{n}'\}$ and $\{\overline{\delta}_{n}\}$ are increasing in $n$ and
$$0< \delta_{n}^{-1}\leq \lambda_{p}\leq \overline{\delta}_{n+1}^{-1}\leq \delta_{n}'^{-1}\leq \sigma_{p}^{-1}< \infty,\qqd  \qd n\geq 1.$$
Furthermore,
$$\lim_{n\rightarrow \infty}\delta_{n}^{-1}=\lim_{n\rightarrow \infty}\overline{\delta}_{n}^{-1}=\lim_{n\rightarrow \infty}\delta_{n}'^{-1}=\lambda_{p}.$$
\end{itemize}}\dethm

 The final assertion of Theorem \ref{t-02} is the main addition to \rf{CWZ}{Theorem 2.4}. The next result is the inverse iteration for obtaining the principal eigenpair of $\Omega_p$ on finite space.

\algo\lb{t-03}\qd{\rm(Inverse iteration)}\qd{\cms Assume that $N<\infty.$ Given a positive $v^{(0)}$ such that $\| v^{(0)}\| _{\mu,p}=1$, at the $n$th$(n \geq 1)$ step, solve the nonlinear equation in $w^{(n)}$
\be \bg{cases}   -\Omega_{p}w^{(n)}(k)=\mu_{k}| v^{(n-1)}_{k}| ^{p-2}v^{(n-1)}_{k}, \qqd k \in E,\\
  w^{(n)}_{N+1}=0.
  \end{cases} \lb{02} \de
and define
$$v^{(n)}=\| w^{(n)}\|^{-1}_{\mu,p}w^{(n)},\qqd   z_{n}=D_{p}(v^{(n)}),$$
where $v^{(n)}_{k}$ denotes the $k$th component of the vector $v^{(n)}$, $\|\cdot\|_{\mu,p}$ denotes the $L^{p}(\mu)$-norm.
Then $\{v^{(n)}\}$ converges to the eigenfunction of $\lambda_{p}.$ Next, $\{z_{n}\}$ is decreasing in $n$ and $$\lim_{n\rightarrow\infty}z_{n}=\lambda_{p}.$$
}\dealgo

As far as we know, the inverse iteration method in numerical algebraic computation was extended to nonlinear eigenproblems by Hein and Bühler in \cite{HB} in 2010.
For $p$-Laplacian, the inverse iteration method was introduced in \cite{BEM} in 2009 by Biezuner, Ercole and Martins.
In 2015, the inverse iteration for obtaining $q$-eigenpairs was given in \cite{EG} by Ercole.
In 2016, the inverse iteration was introduced in \cite{HL} by Hynd and Lindgren for $p$-ground states.
Apart from the ordinary $p$-Laplacian operator, we are studying the weighted one,
as shown by \eqref{02} including two measures $\mu$ and $\nu$. It is quite natural and similar to the generalization from the Laplacian operator to the elliptic one.

Here, we mention that this paper is based on the author's master thesis \cite{LYS17} with some improvements. Besides, further work about shifted inverse iteration was presented in Chen \cite{chen18} where the work related to this paper was cited.

\noindent\\[4mm]

\noindent{\bbb 2\quad Examples}\lb{s-02}\\[0.1cm]
To illustrate the convergence of Theorem \ref{t-02} and Algorithm \ref{t-03}, we consider two examples taken from \rf{CWZ}{Example 2.6 and 2.7}. Using the above two different methods to compute the maximal eigenvalue with MATLAB-R2013a. Here, unless otherwise stated, the initial we take is  $w^{(0)}=f^{(1)}=\hat{\nu}[\cdot,N]^{1/p^{*}}$ as introduced in \rf{CWZ}{\S  2} and $v^{(0)}=\|w^{(0)}\|^{-1} _{\mu,p}w^{(0)}$.

\xmp\qd{\rm\rf{CWZ}{Examples 2.6}}\lb{t-04}\;\;{\rm
Assume that $E=\{0,1,\ldots N\}.$ Let $\mu_{k}=20^{k},$ $\nu_{k}=20^{k+1}$ for $k \in E.$
Then $\lambda_{p}\approx 0.782379$ when  $p=4.5$ and $N=80$. Partial results with the inverse iteration and approximation procedure in Theorem \ref{t-02} are given in Tables 1 and 2. Here and in what follows, we stop our computation once the six precisely significant digits is achieved.}\dexmp

Table 1 and Figure 1 show that $\{z_{k}\}$ decreases faster for smaller $k,$ then goes down slowly and finally becomes stable, which is similar to the linear case in \cite{C1}.

From Tables 1 and 2, it follows that $1/\delta_{n}\leq 1/\overline{\delta}_{n+1}\leq 1/\delta'_{n}$. Besides, it also shows that the sequence $\{z_{n}=1/\overline{\delta}_{n+1}\}$ converges faster than the sequences $\{1/\delta_{n}\}$ and $\{1/\delta'_{n}\}.$

\vspace{0.4truecm}
\begin{center}Table 1. Partial outputs of $z_{k}=\overline{\delta}_{k+1}^{-1}$\end{center}
\vspace{-0.4truecm}
$$\begin{tabular}{|c|c|c|c|c|c|}
\hline
$k$&$z_{k}$&$k$ &$z_{k}$&$k$&$ z_{k}$\\
\hline
0&0.828685&11&0.785009&90&0.782417\\
1&0.796974&12&0.784798&100&0.782403\\
2&0.791961&13&0.784613&120&0.782388\\
3&0.789715&14&0.784449&140&0.782383\\
4&0.788379&15&0.784303&160&0.782381\\
5&0.787472&16&0.784172&180&0.782380\\
6&0.786805&20&0.783759&200&0.782380\\
7&0.786291&30&0.783155&210&0.782380\\
8&0.785878&40&0.782839&$\geq$217&0.782379\\
9&0.785539&50&0.782656&${ }$&${ }$\\
10&0.785253&70&0.782482&${ }$&${ }$\\
\hline
\end{tabular}$$
\begin{center}Table 2. Partial outputs of ($n,$ $1/\delta_{n+1}$ $1/\delta'_{n+1}$)\end{center}
\vspace{-0.6truecm}
$$\begin{tabular}{|c|c|c|c|c|c|}
\hline $n$&$1/\delta_{n+1}$&$1/\delta'_{n+1}$&$n$&$1/\delta_{n+1}$&$1/\delta'_{n+1}$\\
\hline
0&0.778535&4.44787&16&0.778549&0.818055\\
1&0.778535&1.59105&20&0.778579&0.809951\\
2&0.778535&1.14480&30&0.778796&0.799288\\
3&0.778535&1.01136&40&0.779214&0.794022\\
4&0.778535&0.948633&50&0.779722&0.790885\\
5&0.778535&0.912311&70&0.780660&0.787332\\
6&0.778535&0.888671&90&0.781319&0.785410\\
7&0.778535&0.872080&100&0.781552&0.784767\\
8&0.778535&0.859805&120&0.781877&0.783870\\
9&0.778536&0.850360&140&0.782074&0.783315\\
10&0.778536&0.842871&160&0.782194&0.782970\\
11&0.778537&0.836789&180&0.782266&0.782755\\
12&0.778538&0.831752&200&0.782309&0.782621\\
13&0.778539&0.827513&210&0.782324&0.782574\\
14&0.778542&0.823897&217&0.782332&0.782548\\
15&0.778545&0.820776&218&0.782333&0.782544\\
\hline
\end{tabular}$$
\begin{figure}[!h]
\begin{center}{\includegraphics[width=9.0cm,height=5.5cm]{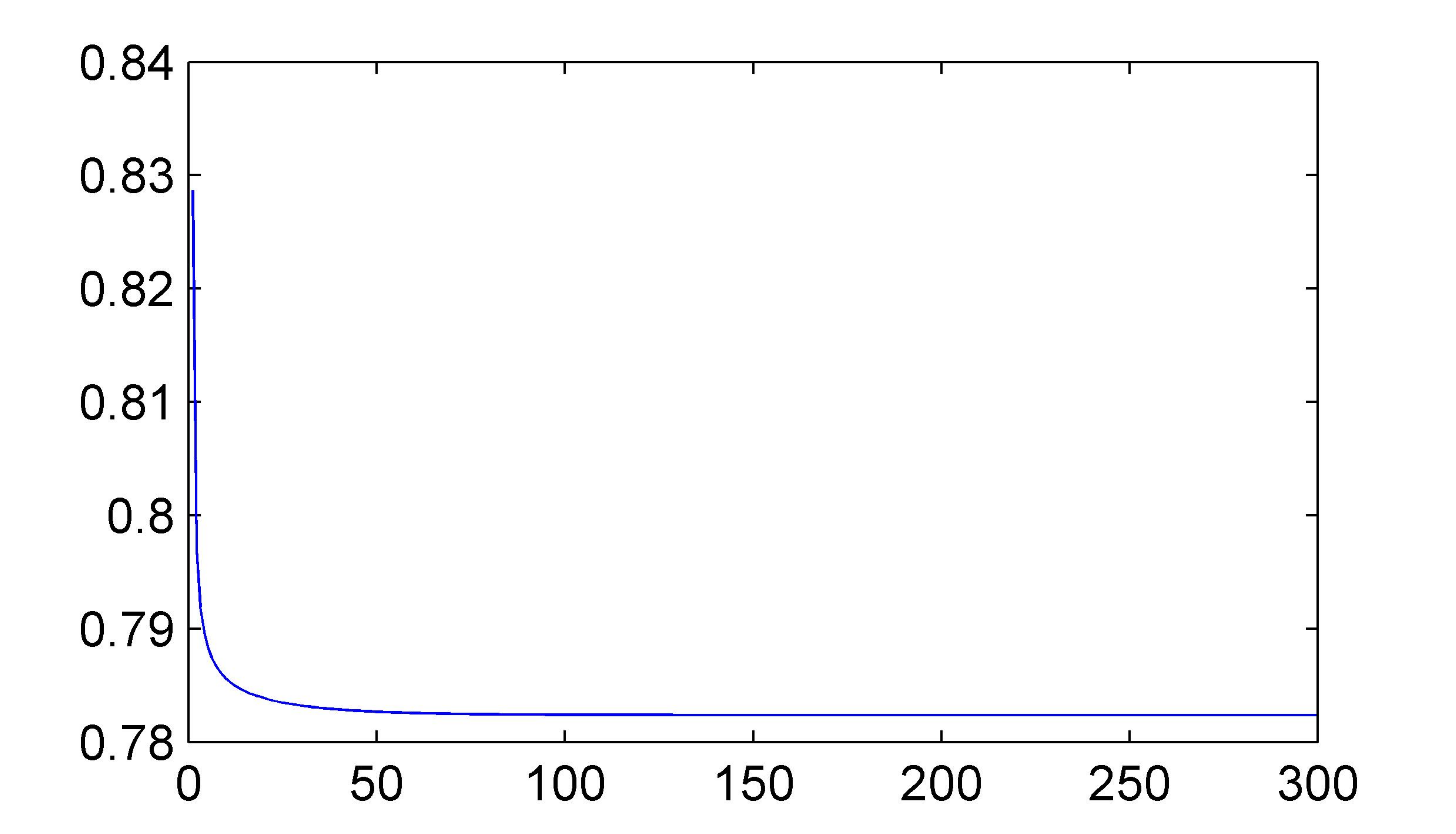}}\end{center}
\vspace{-1.0truecm}
 \begin{center}{\bf Figure 1}\quad The figure of $z_{k}$ for $k=0,1,\ldots, 350$.\end{center}
\end{figure}

To see the serious influence of different initials for inverse iteration, we simply take
 \bgar
w^{(0)} &=& \hat{\nu}[\cdot,N]^{1/p^{*}},\lb{03}\\
w^{(0)} &=& 1,\lb{04}
 \dear
as initials, respectively. The comparing results are given in Table 3.
\vspace{-0.2truecm}
\begin{center}Table 3. Outputs when $N=20$, $p=3$ for different initials\end{center}
\vspace{-0.4truecm}
$$\begin{tabular}{|c|c|c|c|c|c|c|}
\hline
$\displaystyle w^{(0)}$&${\displaystyle z_{0}}$&$z_{1}$&$z_{2}$&$z_{3}$&$z_{4}$&$z_{5}$\\
\hline
$(\ref{03})$
&5.54656&5.30688&5.24358&5.21557&5.19956&5.18922\\

$(\ref{04})$
& 19.  &14.2599&9.95816&7.96854&7.06961&6.58198\\
\hline
$\displaystyle w^{(0)}$&$z_{6}$&$z_{7}$&$z_{8}$&$z_{9}$&$z_{10}$&$z_{15}$\\
\hline
 $(\ref{03})$
&5.18209&5.17696&5.17318&5.17036&5.16822&5.16311\\

 $(\ref{04})$
&6.27831&6.07182&5.92271& 5.81020&5.72241&5.47138\\
\hline
$\displaystyle w^{(0)}$&$z_{20}$&$z_{25}$&$z_{30}$&$z_{40}$&$z_{50}$&$z_{59}$\\
\hline
 $(\ref{03})$
&5.16175&5.16137&5.16127&5.16123&5.16122&5.16122\\

 $(\ref{04})$
&5.35331&5.28531&5.24194&5.19412&5.17383&5.16636\\
 \hline
$\displaystyle w^{(0)}$&$z_{80}$&$z_{100}$&$z_{140}$&$z_{147}$&{ }&{ }\\
\hline
$(\ref{04})$&5.16182&5.16130&5.16122&5.16122&{ }&{ }\\
\hline
\end{tabular}$$

Table 3 shows that a good initial for the inverse iteration is essential.

\xmp\qd{\rm \rf{CWZ}{Example 2.7}}\lb{t-05}\;\;{\rm
Assume that $E=\{0,1,\cdots,N\}.$ Let $\mu_{k}=\nu_{k}=1$ for each $k\in E$. When $p=2.5$ and $N=40,$ the outputs of $z_{n},$ $1/\delta_{n+1},$ $1/\delta'_{n+1},$ and $1/\overline{\delta}_{n+1}$ are given in Table 6.
\begin{center}Table 6. Outputs of ($n,$ $z_{n},$ $1/\delta_{n+1},$ $1/\delta'_{n+1},$ $1/\overline{\delta}_{n+1}$)\end{center}

$$\begin{tabular}{|c|c|c|c|c|}
\hline
     $\displaystyle n$&$z_{n}$  &$\displaystyle \delta^{-1}_{n+1}$&$\displaystyle \delta'^{-1}_{n+1}$&$\displaystyle \overline{\delta}^{-1}_{n+1}$\\
\hline
      0&0.000458009&	0.000255829&	0.00160159&	0.000458009\\

      1&0.000271491&	0.000269664&	0.000283578&0.000271491\\

      2&0.000271279&	0.000271048&	0.000272059&0.000271279\\

      3&0.000271277&	0.000271250&	0.000271349&0.000271277\\

      4&${  }$&	        0.000271274&	0.000271285&${  }$\\

      5&${  }$&	        0.000271277&	0.000271278&${  }$\\

      6&${  }$&	        ${  }$&	        0.000271277&${  }$\\
\hline
\end{tabular}$$

Here take $z_n$  for instance, when $n\geq 3,$ the outputs are the same and so are omitted.}

\dexmp

\noindent\\[4mm]

\noindent{\bbb 3\quad Proofs of the main results in Section \ref{s-01}}\lb{s-03}\\[0.1cm]
To prove Algorithm \ref{t-03}, we need to solve equation \eqref{02}. Fortunately, the solution can be expressed explicitly in terms of operator $I\!I$, as shown in Lemma \ref{t-06} below. Actually, equation \eqref{02} is a particular case of Poisson Equation \eqref{05}:
  \be
  -\Omega_{p}g(k)=\mu_{k}|f_{k}|^{p-1}\text{sgn}(f_{k}),\qqd k\in E,\lb{05}
  \de
where $ g,f \in \{g:E\rightarrow \mathbb{R} | g_{-1}=g_{0},g_{N+1}=0\}.$ In this section, we assume that $N<\infty.$

\lmm\lb{t-06}\qd{\cms  Given $f>0,$
then the unique solution to \eqref{05} is

\begin{align*}
g(k)=f_k\left(I\!I_{k}(f)\right)^{p^{*}-1},\qqd k\in E.
\end{align*}
}\delmm
\prf
With $f>0$ at hand, by \eqref{01}, making a summation on both sides of equation \eqref{05} over $0\sim k$, we get
  \be
 -\nu_{k}|\partial_{k}(g)|^{p-1}\text{sgn}(\partial_{k}(g))=
 \sum_{j=0}^{k}\mu_{j}f_{j}^{p-1}>0, \qqd k\in E.\lb{06}
 \de
Thus $\partial_{k}(g)<0.$ Combining with the boundary condition $g_{N+1}=0,$ we have $g_{k}>0$ for each $k\in E.$ Besides, by \eqref{06} and the boundary condition at $N$, we have
$g_{N}=\hat{\nu}_{N}\left(\sum_{j=0}^{N}\mu_{j}f_{j}^{p-1}\right)^{p^{*}-1}.$
Now with $g>0$ and $\partial(g)<0$ at hand, by \eqref{06} and induction on $k,$ it follows that
$$g_{k}=\sum_{i=k}^{N}\hat{\nu}_{i}\left(\sum_{j=0}^{i}\mu_{j}f_{j}^{p-1}\right)^{p^{*}-1},\qqd k\in E.$$
Therefore, we conclude $g(k)=f_k\left(I\!I_{k}(f)\right)^{p^{*}-1}$, it is easy to verify the uniqueness of the solution, which proves the lemma.
\deprf

Now in our setup, $\lambda_{p}> 0$ by {\rm \rf{CWZ}{Theorem 2.3}},
and the eigenfunction corresponding to $\lambda_{p}$ is strictly monotone by {\rm \rf{CWZ}{Proposition 3.2}}.
Note that in Algorithm \ref{t-03}, we choose the positive initial $v^{(0)}.$ Now, by induction on $n,$
both $\{w^{(n)}\}$ and $\{v^{(n)}\}$ are positive and decreasing functions according to Lemma \ref{t-06}.
Hence, as the mimic eigenfunctions, we are happy to get the monotone sequences $\{w^{(n)}\}$ and $\{v^{(n)}\}.$
Therefore, in what follows, the modulus can be discarded. Besides, noting that $I\!I(f)=I\!I(cf)$ for each constant $c>0,$ we have Remark \ref{t-07} by induction on $n$.

\rmk\lb{t-07}\qd{\cms Assume that $v^{(0)}=c f^{(1)}$ for a constant $c>0$, where $v^{(0)}$ and $f^{(1)}$ are given in Algorithm \ref{t-03} and Theorem \ref{t-02}, respectively. Then
$$z_n=\bar{\dz}_{n+1}^{-1},\qqd n\geq 0.$$}\dermk
Noticing that $\{v^{(n)}\}$ in Algorithm $\ref{t-03}$ is an approximation sequence of the eigenfunction of $\lambda_{p}$, we have for large enough $n$,
$$-\Omega_{p}v^{(n)}\approx \lambda_{p}\text{diag}(\mu)(v^{(n)})^{p-1}.$$
Furthermore,
$$\lambda_{p}\approx \frac{(-\Omega_{p}v^{(n)},v^{(n)})}{(\text{diag}(\mu)(v^{(n)})^{p-1},v^{(n)})}=\frac{D_{p}(v^{(n)})}{\| v^{(n)}\| _{\mu,p}^{p}}=\overline{\delta}_{n}^{-1}.$$
Hence the convergence of Algorithm $\ref{t-03}$ could be expected.

Define
\be\xi_{n-1}=\frac{1}{\| {w}^{(n)}\|^{p-1}_{\mu,p}}, \qqd n\geq 1.\lb{07}\de
Recall that $v^{(n)}=\|w^{(n)}\|_{\mu,p}^{-1}w^{(n)}.$ Inserting this into \eqref{02}, it follows that
\be
 \begin{cases}
  -\Omega_{p}v^{(n)}(k)=\xi_{n-1}\mu_{k}| v^{(n-1)}_{k}| ^{p-2}v^{(n-1)}_{k},
  \qqd k \in E,\\
 v^{(n)}_{N+1}=0,
  \end{cases}\lb{08}
\de
where $v^{(n)}$ and $v^{(n-1)}$ are given in Algorithm \rm{\ref{t-03}}.

Now we start to prove Algorithm \ref{t-03}. Inspired by the method used in {\rf{HL}{Proposition 2.4}}, we have the following proposition for the convergence of $\{z_{n}\}.$
\prp\lb{t-08}\qd{\cms The sequences
$$ z_{n}=D_{p}(v^{(n)})
\qqd  \text{\cms and}\qqd \xi_{n}=\| {w}^{(n+1)}\|^{1-p}_{\mu,p} $$
defined in Algorithm \ref{t-03} and \eqref{07}, respectively, are decreasing in $n$
and satisfy
$$\lambda_{p}\leq z_{n+1}\leq \xi_{n} \leq z_{n} \leq \xi_{n-1}.$$
Therefore
\be\lz_p\leq\lim_{n\rightarrow \infty}z_{n}=\lim_{n\rightarrow\infty }\xi_{n}=:\xi. \lb{09}\de
}\deprp

\prf Since $\nu_{-1}=0$ and $f_{N+1}=0,$ we have for arbitrary $\{H_{k}\}$
$$\sum_{k=0}^{N}\nu_{k-1}f_{k}H_{k-1}=\sum_{k=-1}^{N-1}\nu_{k}f_{k+1}H_{k}
=\sum_{k=0}^{N}\nu_{k}f_{k+1}H_{k}.$$

Combining this with \eqref{01}, we have
\begin{align}
    \quad& (-\Omega_{p}g,f)\nonumber\\
     &=-\sum_{k=0}^{N}\nu_{k}f_{k}| \partial_{k}(g)|^{p-1}\text{sgn}\left(\partial_{k}(g)\right)+\sum_{k=0}^{N}\nu_{k-1}f_{k}| \partial_{k-1}(g)|^{p-1}\text{sgn}\left(\partial_{k-1}(g)\right) \qd \nonumber\\
     &=\sum_{k=0}^{N}\nu_{k}| \partial_{k}(g)|^{p-1}\text{sgn}\left(\partial_{k}(g)\right)\partial_{k}(f). \lb{10}
   \end{align}
In particular, $$(-\Omega_{p}g,g)=\sum_{k=0}^{N}\nu_{k}| \partial_{k}(g)|^{p}=D_{p}(g).$$
If $(g,f)$ satisfy Poisson equation \eqref{05} and $f\neq 0$,
then for each $h \in \{g:E\rightarrow \mathbb{R} | g_{-1}=g_{0},g_{N+1}=0\},$ we have
\begin{align*}
 |(-\Omega_{p}g,h)|
 \overset{\eqref{05}}{=}|\left(|f|^{p-1}\text{sgn}(f),h\right)_{\mu}|
 \leq\| f\|_{\mu,p}^{p-1} \| h \|_{\mu,p}
 \qd(\text{H\"{o}lder inequality}).
\end{align*}
Dividing both sides by $\|g\|_{\mu,p}^{p},$ it follows that
$$\frac{|(-\Omega_{p}g,h)|}{\|g\|_{\mu,p}^{p}}\leq
\frac{\| f\|_{\mu,p}^{p-1} \| h \|_{\mu,p}}{\|g\|_{\mu,p}^{p}}.$$
In particular, when $h=g$, we have
\be
\frac{D_{p}(g)}{\| g\|_{\mu,p}^{p} }\leq\frac{\| f\|_{\mu,p}^{p-1} }{\|g\|_{\mu,p}^{p-1} }
=\frac{\|f\|_{\mu,p}^{p} }{\| g\|_{\mu,p}^{p-1} \| f\|_{\mu,p}}. \lb{11}
\de
Note that here we use the norm $\|\cdot\|_{\mu,p}$ to control the ordinary inner product (in particular, control $D_{p}$). Now, we use $D_{p}$ to control the norm $\|\cdot \|_{\mu,p}$. Again, let $(g,f)$ satisfy Poisson equation \eqref{05} and $f\neq 0$, then
\be
 \begin{split}
\| f\|_{\mu,p}^{p}= \left(|f|^{p-1}\text{sgn}(f) ,f\right)_{\mu}&\overset{\eqref{05}}{=}(-\Omega_{p}g,f)\qqd \\
&\overset{\eqref{10}}{=}\sum_{k=0}^{N}\nu_{k}|\partial_{k}(g)|^{p-1}\text{sgn}\left(\partial_{k}(g)\right)\partial_{k}(f)\\
&=\left(|\partial(g)|^{p-1}\text{sgn}\left(\partial(g)\right),\partial(f)\right)_{\nu}\\
&\overset{\text{H\"{o}lder~inequality}}{\leq} D_{p}(g)^{\frac{p-1}{p}} D_{p}(f)^{\frac{1}{p}}.\lb{12}
 \end{split}
\de
Combining \eqref{11} with \eqref{12}, it follows that
$$\frac{D_{p}(g)}{\| g\|_{\mu,p}^{p} }\leq\frac{D_{p}(g)^{\frac{p-1}{p}} D_{p}(f)^{\frac{1}{p}}}{\| g\|_{\mu,p}^{p-1}\| f\|_{\mu,p}}.$$
By a cancellation procedure, we get
$$
\frac{D_{p}(g)}{\| g\|_{\mu,p}^{p} }\leq\frac{D_{p}(f)}{\| f\|_{\mu,p}^{p} }.
$$
Now letting  $f=v^{(n)},~ g=w^{(n+1)}$ and noting that $\| v^{(n)}\|_{\mu,p}=1,$ we have
\begin{align*}
z_{n+1}=\frac{D_{p}(w^{(n+1)})}{\| w^{(n+1)}\|_{\mu,p}^{p} } &\overset{\eqref{11}}{\leq}
\frac{1}{\| w^{(n+1)}\|_{\mu,p}^{p-1}}\frac{\| v^{(n)}\|_{\mu,p}^{p} }{\| v^{(n)}\|_{\mu,p}}
=\frac{1}{\| w^{(n+1)}\|_{\mu,p}^{p-1}}=\xi_{n}\\
&\overset{\eqref{12}}{\leq}
\frac{D_{p}(w^{(n+1)})^{\frac{p-1}{p}} D_{p}(v^{(n)})^{\frac{1}{p}}}{\| w^{(n+1)}\|_{\mu,p}^{p-1}\| v^{(n)}\|_{\mu,p}}=z_{n+1}^{(p-1)/p}z_{n}^{1/p},
\end{align*}
By a simple cancellation procedure, we get
$$z_{n+1}\leq \xi_{n}\leq z_{n}\leq \xi_{n-1}.$$
This implies not only the decreasing property of $\{z_n\}$ and $\{\xi_n\}$ in Algorithm \ref{t-03} and \eqref{07}, respectively, but also \eqref{09}.
\deprf

Next, we prove the convergence of the mimic eigenfunction sequence $\{v^{(n)}\}$.

\prp\lb{t-09}\qd{\cms The sequence $\{v^{(n)}\}$ converges to an eigenfunction of $\Omega_p$ corresponding to $\xi.$
}\deprp

\prf To prove the result, we adopt two steps.
First, we prove that there exists a subsequence of $\{v^{(n)}\}$ converging to an eigenfunction of $\Omega_p$ corresponding to $\xi.$
Next, we prove that any convergent subsequence of $\{v^{(n)}\}$ converges to the same $v.$

(a) Prove that there exists a subsequence of $\{v^{(n)}\}$ converging to an eigenfunction of $\xi.$

Since $\{v^{(n)}\}$ is on the unit sphere of $L^{p}(\mu)$,
there exists a subsequence $\{v^{(n_{k})}\}$ of $\{v^{(n)}\}$ and a function $v$ on $E$ such that
$$v^{(n_{k})}\longrightarrow v \qqd \text{(pointwise)}.$$
Thus it suffices to prove that $v$ is an eigenfunction of $\xi.$ By Lemma \ref{t-06}, $v^{(n)}~(n\geq 1)$ is positive and strictly decreasing, therefore, $v$ is nonnegative and decreasing.
According to \eqref{08}, we get
$$-\Omega_{p}v^{(n_{k}+1)}(i)=
\xi_{n_{k}}\mu_{i} \left(v^{(n_{k})}_{i}\right)^{p-1}, \qquad i \in E.$$
Again, by Lemma \ref{t-06} $$v_{\ell}^{(n_{k}+1)}=\sum_{i=\ell}^{N}\hat{\nu}_{i}
\left(\sum_{j=0}^{i}\xi_{n_{k}}\mu_{j}\left(v_{j}^{(n_{k})}\right)^{p-1}\right)^{p^{*}-1}.$$
By Proposition \ref{t-08}, we have $\xi_{n_{k}}\rightarrow \xi$. Combing this with $v^{(n_{k})}\rightarrow v \text{(pointwise)},$
it follows that
$$v^{(n_{k}+1)}\longrightarrow \bar{v} \qqd  \text{(pointwise)},$$
where
\be
\bar{v}_{\ell}=\sum_{i=\ell}^{N}\hat{\nu}_{i}\left(\sum_{j=0}^{i}\xi\mu_{j}v_{j}^{p-1}\right)^{p^{*}-1}.\lb{13}
\de
Since
$$\lim_{k\rightarrow\infty}z_{n_{k}}=\lim_{k\rightarrow\infty}D_{p}\left(v^{(n_{k})}\right)=\xi
=\lim_{k\rightarrow\infty}z_{n_{k}+1}=\lim_{k\rightarrow\infty}D_{p}\left(v^{(n_{k}+1)}\right),$$
and $$\|v^{(n_{k})}\|_{\mu,p}=\|v^{(n_{k}+1)}\|_{\mu,p}=1,$$
we get
$D_{p}(\bar{v})=D_{p}(v)=\xi, \|\bar{v}\|_{\mu,p}=\|v\|_{\mu,p}=1$.
Combining with \eqref{13}, it follows that $\bar{v}$ is positive and decreasing.
 Besides, it can easily get that
\be
-\Omega_{p}\bar{v}(\ell)=\xi\mu_{\ell}v_{\ell}^{p-1}.\lb{014}
\de
Furthermore, we have
$$\xi =\xi \|v\|_{\mu,p}^{p}=(\xi v^{p-1},v)_{\mu}\overset{\eqref{014}}{=}(-\Omega_{p}\bar{v},v)
\overset{\eqref{10}}{=}(|\partial(\bar{v})|^{p-1}\text{sgn}(\partial{\bar{v}}),\partial v)_{\nu},$$
by H\"{o}lder~inequality, we have
\begin{equation}
(|\partial(\bar{v})|^{p-1}\text{sgn}(\partial{\bar{v}}),\partial v)_{\nu}
\leq D_{p}(\bar{v})^{\frac{p-1}{p}} D_{p}(v)^{\frac{1}{p}}=\xi,\lb{15}
\end{equation}
which means the equality in \eqref{15} holds. Thus

there exist constants $c_{1}$ and $c_{2}$ satisfying $c_{1}c_{2}\neq 0$ such that
~$c_{1}|\partial{\bar{v}}|^{(p-1)\cdot p/(p-1)}+c_{2}|\partial{v}|^{p}=0.$
~Combining this with
$D_{p}(v)=D_{p}(\bar{v})=\xi,$
~we get~$|\partial_{j}{\bar{v}}|=|\partial_{j}v|.$
Since both $v$ and $\bar{v}$ are decreasing, it follows that
~$\partial_{j}{\bar{v}}=\partial_{j}v$.
By the boundary condition at $N+1$, we get
~$\bar{v}_N=v_N.$
With positive $\bar{v}$ at hand, by induction on $k$, we conclude that
$$v_k=\bar{v}_k> 0,\qqd k\in E.$$
Combing the positivity and decreasing property of $\bar{v}$ with \eqref{014}, it  follows that $v$ is a positive and decreasing eigenfunction of $\xi$. Hence
$$-\Omega_{p}v{(k)}=\xi\mu_{k}v_{k}^{p-1}.$$

(b) Prove that any convergent subsequence of $\{v^{(n)}\}$ converges to the same $v.$

In fact, the proof in (a) means that $v$ is a solution to
\be\begin{cases}
 -\Omega_{p}f=\xi\text{diag}(\mu) f^{p-1},\\
 f_{N+1}=0.
 \end{cases}\lb{16}\de
It is easy to verify that equation \eqref{16} has a unique positive and decreasing solution $f$ with $\|f\|_{\mu,p}=1.$
Thus every convergent subsequence of $\{v^{(n)}\}$ converges to the same $v$.

The proofs (a) and (b) above show that the sequence $\{v^{(n)}\}$ converges to an eigenfunction of $\xi$ and hence $(\xi,v)$ is an eigenpair of $-\Omega_{p}.$
\deprf

By Propositions \ref{t-08} and \ref{t-09}, the sequence $\{(v^{(n)}, z_n)\}$ is an approximation one for eigenpair of $-\Omega_p$. As the mimic eigenfunctions, $\{v^{(n)}\}$ is a positive and decreasing sequence which has the same property with the principal eigenfunction of $-\Omega_p.$ Thus, it is expected to prove that $\xi=\lz_p.$ In fact, the double summation form of variational formulas in \cite{CWZ} plays a typical role to prove that $\{(v^{(n)}, z_n)\}$ is an approximation sequence for the principal eigenpair of $-\Omega_p$.
\prp\lb{t-10}\qd{\cms Each of the sequences $\{\delta_{n}\},\{\delta_{n}'\}$ and $\{\bar{\delta}_{n}\}$ in Theorem \ref{t-02} converges to $\lambda_{p}$,
$$\lim_{n\rightarrow\infty}\delta_{n}^{-1}=\lim_{n\rightarrow \infty}\delta_{n}'^{-1}
=\lim_{n\rightarrow\infty}\bar{\delta}_{n}^{-1}=\lambda_{p}=\xi.$$
}\deprp
\prf
(a)~First, we prove $$\lim_{n\rightarrow\infty}\frac{w^{(n)}}{w^{(n+1)}}=1\qqd \text{(pointwise)}.$$

By Proposition \ref{t-09}, $\{v^{(n)}\}$ converges to an eigenfunction (denoted by $v$) of $\xi$  and $v$ is positive and decreasing. Thus by Proposition \ref{t-08}, we have
 $$\lim_{n\rightarrow\infty}\frac{w^{(n)}}{w^{(n+1)}}=\lim_{n\rightarrow\infty}\frac{v^{(n)}}{v^{(n+1)}}
\left(\frac{\| w^{(n)}\|_{\mu,p}^{1-p}}{\| w^{(n+1)}\|_{\mu,p}^{1-p}}\right)^{1-p^*}\!\!\!\!= 1\qd \qd \text{(pointwise)}.$$

(b) Prove $$\lim_{n\rightarrow\infty}I\!I(w^{(n)})=\lim_{n\rightarrow\infty}I\!I(v^{(n)})=\xi^{-1}.$$

In fact, since
$$w^{(n+1)}=v^{(n)}\left(I\!I(v^{(n)})\right)^{p^{*}-1}=\|w^{(n)}\|_{\mu,p}^{-1}w^{(n)}
\left(I\!I(v^{(n)})\right)^{p^{*}-1},$$
we get
$$\lim_{n\rightarrow\infty}I\!I(v^{(n)})=\lim_{n\rightarrow\infty}\bigg(w^{(n+1)}/w^{(n)}\bigg)^{p-1}\|w^{(n)}\|_{\mu,p}^{p-1}
=\xi^{-1}.$$

(c) Prove $$\xi=\lambda_{p}.$$

In fact, (b) means that $$\xi^{-1}=\lim_{n\rightarrow\infty}I\!I(v^{{(n)}})=\lim_{n\rightarrow \infty}\sup_{k \in E}I\!I_{k}(v^{(n-1)}).$$

Since $I\!I$ is homogeneous with respect to a constant, let $v^{(0)}$ coincide with $f^{(1)}$ up to a fixed constant, then
$$
\lim_{n\rightarrow \infty}\sup_{k \in E}I\!I_{k}(v^{(n-1)})=\lim_{n\rightarrow \infty}\sup_{k \in E}I\!I_{k}(f^{(n)})
=\lim_{n\rightarrow \infty}\delta_{n},$$
which means $\lim_{n\rightarrow \infty}\delta_{n}=\xi^{-1}.$ For similar reason,
$$ \xi^{-1}=\lim_{n\rightarrow \infty}\inf_{k \in E}I\!I_{k}(v^{(n-1)})=
\lim_{n\rightarrow \infty}\inf_{k \in E}I\!I_{k}(f^{(n)})
=\lim_{n\rightarrow \infty}\delta'_{n}.$$
By {\rm \rf{CWZ}{Theorem 2.1}}, we know that
\be\delta_{n}^{-1}\leq \lambda_{p}\leq \delta_{n}'^{-1}.\lb{017}\de
Taking limit on both sides of $(\ref{017})$ and noting that
$$\lim_{n\rightarrow \infty}\delta_{n}^{-1}=\lim_{n\rightarrow \infty}\delta_{n}'^{-1}=\xi,$$
 we get
$\xi=\lambda_{p}.$ Combing with Remark \ref{t-07}, we get $\lim\limits_{n\rightarrow\infty}\bar{\dz}_n=1/\lz_p.$
\deprf

Now, we have proved Algorithm \ref{t-03} in terms of Propositions \ref{t-08} $\sim$\ref{t-10}. For Theorem \ref{t-02}, the remained results can be obtained easily from the dual variation of $\lambda_{p}$ in \cite{CWZ}, refer to {\rm \rf{CWZ}{Theorem 2.1 and Theorem 2.4}}. Hence, we have also completed the proof of Theorem \ref{t-02}. Here, we mention that in Theorem \ref{t-02} and Algorithm \ref{t-03}, we need only $v^{(0)}>0$, which can be seen from the linear case according to Perron-Frobenius theorem. The linear case can be proved beautifully using spectrum decomposition.
\begin{figure}[ht]\vspace{0.5cm}
\begin{center}{\includegraphics[width=11.0cm,height=7.5cm]{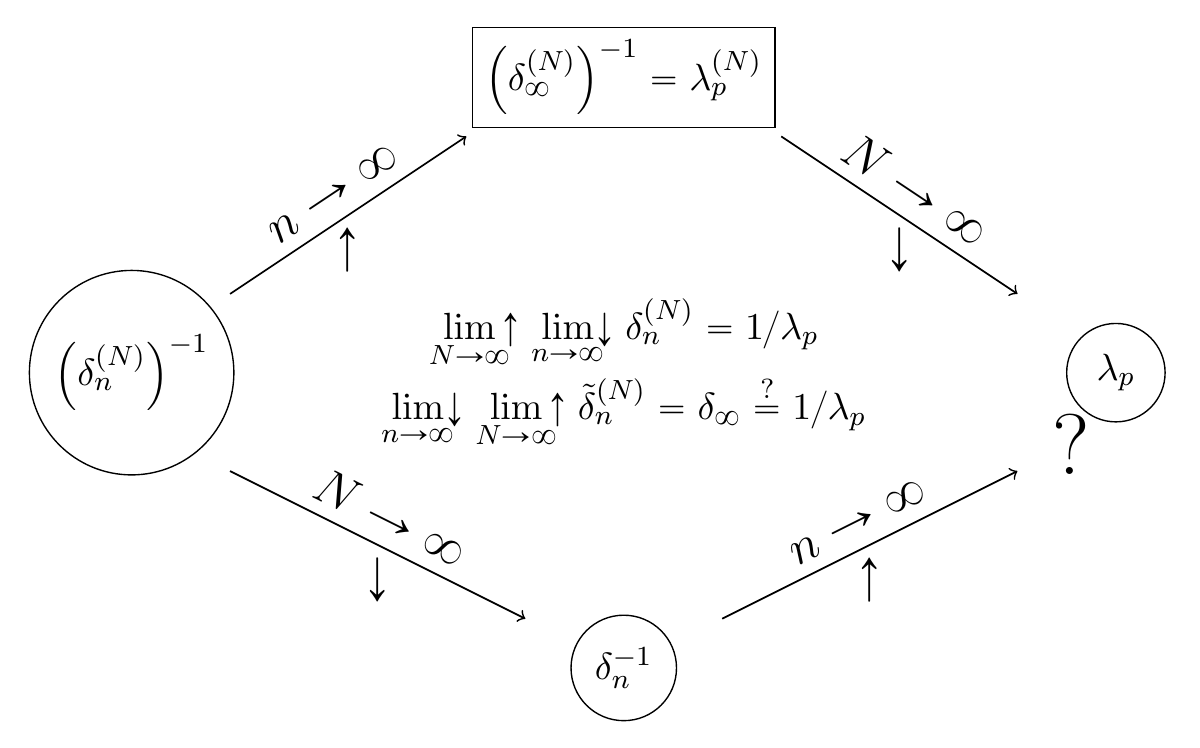}}\end{center}
 \begin{center}{\bf Diagram 1}\quad Here $\uparrow~(\downarrow)$ denotes increasing~(decreasing).\end{center}
\end{figure}

To conclude this section, we discuss the infinite case: $N=\infty$.
Roughly speaking, what we know is shown by Diagram 1. To explain the meaning of the
diagram, noticing 
\bg{itemize}\setlength{\itemsep}{-0.6ex}
\item[(a)] by {\rm \rf{CWZ}{Theorem 2.4(1)}}, for fixed $N,$ $\{\delta_n^{(N)}\}$
is decreasing in $n$; 
\item[(b)] by Proposition \ref{t-10}, $\lim\limits_{n\rightarrow\infty}\delta_n^{(N)}=1/\lz^{(N)}_p$;
\item[(c)] by {\rm \rf{CWZ}{Remark 3.3(2)}}, $\lim\limits_{N\rightarrow\infty}\lz_p^{(N)}=\lz_p.$
\end{itemize}
Combining these three facts together, we obtain
$$\lim\limits_{N\rightarrow\infty}\!\!\!\uparrow\lim\limits_{n\rightarrow\infty}\!\!\!\downarrow \delta_{n}^{(N)}=1/\lz_p.$$

On the other hand, by Theorem \ref{t-01}\,(1), we have a sequence $\{f^{(n)}\}_{n\geq 1}$ on $E=\{k\in \mathbb{Z}_+: 0\leq k<\infty\}$. Hence, we can define
$${\tilde \delta}_n^{(N)}=\sup_{0\leq i\leq N}I\!I_i(f^{(n)}\mathbb{1}_{\leq N}).$$
Then it is easy to check that $\{{\tilde\delta}_n^{(N)}\}$ is increasing in $N,$ and its limit is $\dz_n$ defined in Theorem \ref{t-01}\,(1) when $N=\infty$. Furthermore, by Theorem \ref{t-01}\,(1), $\{\dz_n\}$ is decreasing in $n,$ denoting its limits by $\dz_{\infty}$:
$$\dz_{\infty}=\lim\limits_{n\rightarrow\infty}\!\!\!\downarrow
\lim\limits_{N\rightarrow\infty}\!\!\!\uparrow {\tilde \delta}_{n}^{(N)}.$$
Therefore, the main open question here is $\dz_{\infty}=1/\lz_p$? It is regretted that
at the moment we do not know how to solve this problem.

For the dual boundary, similar results hold, which are stated in the next section.
\noindent\\[4mm]

\noindent{\bbb 4\quad Results on the dual boundary condition(DN-case)}\lb{s-04}\\[0.1cm]
In this section, the notation are the same as those used in \rf{CWZ}{\S  4}, unless otherwise stated.
Let $E:=\{k\in \mathbb{Z}_+: 0\leq k < N+1\}$ with $N\leq\infty,$  $\partial^{-}_{k}(f)=f_{k-1}-f_{k}$, where $\partial^{-}_{k}$ denotes the left difference, then $\Omega_{p}$ can be rewritten as
$$
\Omega_{p}f(k)=-\nu_{k+1}|\partial_{k+1}^{-}(f)|^{p-1}\text{sgn}(\partial^{-}_{k+1}(f))+\nu_{k}
|\partial^{-}_{k}(f)|^{p-1}\text{sgn}(\partial^{-}_{k}(f)).
$$
Under the DN-boundary condition, the operator $I\!I$ is defined as follows
\begin{align*}
I\!I_{i}(f)=\frac{1}{f_{i}^{p-1}}\left[\sum_{j=1}^{i}\hat{\nu}_{j}
\left(\sum_{k=j}^{N}\mu_{k}f_{k}^{p-1}\right)^{p^{*}-1}\right]^{p-1},\qqd i\in E,
\end{align*}
and $\sigma_{p}$ is defined below $$\sigma_{p}=\sup_{n \in E }\left(\mu[n,N]\hat{\nu}[1,n]^{p-1}\right).$$
Here we mention that the positive sequence $\{\nu_k:k\in E\}$ satisfies the boundary condition that $\nu_{N+1}:=0$ if $N<\infty$. When $N\leq\infty,$ the approximation procedure for $\lz_p$ is presented in Theorem \ref{t-11}.

\thm\lb{t-11}\qd{\cms Assume that $N\leq\infty$  and $\sigma_{p}<\infty$.
\bg{itemize} \setlength{\itemsep}{-0.6ex}
\item[(1)] Define
$$f^{(1)}=\hat{\nu}[1,\cdot]^{1/p^*},\qd f^{(n)}=f^{(n-1)}\left(I\!I\left(f^{(n-1)}\right)\right)^{p^{*}-1},\qd (n\geq 2),$$
and $\dz_n=\sup\limits_{i\in E}I\!I_i\left(f^{(n)}\right).$ Then $\dz_n$ is decreasing in $n$~(denote its limit by $\dz_{\infty}$) and
$$\lz_p\geq \dz_{\infty}^{-1}\geq\cdots\geq\dz_1^{-1}>0.$$
\item[(2)] For fixed $m \in E$, define
$$f^{(1,m)}\!\!=\!\!\hat{\nu}[1,\cdot\wedge m],\qd
f^{(n,m)}\!\!=\!\!f^{(n-1,m)}\left(I\!I\left( f^{(n-1,m)}\right)(\cdot\wedge m)\right)^{p^*-1},\qd n\geq 2,$$
and
$$\dz_n^{'}=\sup_{m\in E}\inf_{i\in E}I\!I_i\left(f^{(n,m)}\right).$$
Then $\dz_n^{'}$ is increasing in $n$~(denote its limit by $\dz_{\infty}^{'}$)~and
$$\sigma_{p}^{-1}\geq\dz_1^{'-1}\geq\cdots\geq\dz_{\infty}^{'-1}\geq\lz_p.$$
Next, define
$$\bar{\dz}_n=\sup_{m\in E}\frac{\mu\left(f^{(n,m)^p}\right)}
{D_{p}\left(f^{(n,m)}\right)},\qqd n\geq 1.$$
Then $\bar{\dz}_n^{-1}\geq\lz_p$ and $\bar{\dz}_{n+1}\geq\dz_{n}^{'}$ for every $n\geq 1.$
\end{itemize}}\dethm

Theorem \ref{t-11} was first presented in \rf{CWZ}{Theorem  4.3}. There also exist three similar open questions in Theorem \ref{t-11} and the answers for finite $N$ are stated in Theorem \ref{t-12} below.
\thm\lb{t-12}\qd{\cms Assume that $N<\infty.$ Let $f^{(1)}>0$ be given.
\ben \setlength{\itemsep}{-0.6ex}
\item[(1)] Define successively functions on $E:$

$$f^{(n)}=f^{(n-1)}\left(I\!I\left(f^{(n-1)}\right)\right)^{p^{*}-1},\qqd n\geq 2.$$
\item[(2)] For $n\geq 1$, define three sequences as follows
 $$\delta_{n}=\sup_{i \in E} I\!I_{i}\left(f^{(n)}\right),\qqd \delta_{n}'=\inf_{i\in E}I\!I_{i}\left(f^{(n)}\right),\qqd \bar{\delta}_{n}=\frac{\mu({f^{(n)^{p}}})}{ D_{p}({f^{(n)}})}.$$
Then $\{\delta_{n}\}$ is decreasing in $n$, both $\{\delta_{n}'\}$ and $\{\bar{\delta}_{n}\}$ are increasing in $n$ and
$$0< \delta_{n}^{-1}\leq \lambda_{p}\leq \bar{\delta}_{n+1}^{-1}\leq \delta_{n}'^{-1}\leq \sigma_{p}^{-1}< \infty,\qqd \text{\rm for any} \qd n\geq 1.$$
Furthermore,
$$\lim_{n\rightarrow \infty}\delta_{n}^{-1}=\lim_{n\rightarrow \infty}\bar{\delta}_{n}^{-1}=\lim_{n\rightarrow \infty}\delta_{n}'^{-1}=\lambda_{p}.$$
\den
}\dethm

The final assertion of Theorem \ref{t-12} is the main addition to \rf{CWZ}{Theorem  4.3}.
The next result is the inverse iteration for obtaining the principal eigenpair of $\Omega_p$ corresponding to DN-boundary condition.

\algo\lb{t-13}\qd{\rm(Inverse iteration)}\qd \cms{Assume that $N<\infty.$ Given a positive $w^{(0)}$, at the $(n+1)$th step, let $w^{(n+1)}$ be the solution to the following equation
\be
 \begin{cases}
  -\Omega_{p}w(k)=\mu_{k}| w^{(n)}_{k}| ^{p-2}w^{(n)}_{k}, \qqd k \in E.\\
  w_{0}=0.
  \end{cases} \lb{018}
\de
Let$$z_{n+1}=\frac{D_{p}(w^{(n+1)})}{\|w^{(n+1)}\|^{p}_{\mu,p}}.$$
Then the sequence $\{w^{(n)}/\|w^{(n)}\|_{\mu,p}\}$ converges to the eigenfunction of $\lz_{p}.$ Next, $\{z_{n}\}$ is decreasing in $n$ and $$\lim_{n\rightarrow\infty}z_{n}=\lambda_{p}.$$
}\dealgo

Since the proof here is similar to the one given in Section 3, the only different point is that the sum should be the tail summation due to the boundary condition. Let $v^{(n)}=\|w^{(n)}\|^{-1}_{\mu,p}w^{(n)},$ the equivalent form of \rm{\eqref{018}} is
\begin{equation*}
 \begin{cases}
  -\Omega_{p}v^{(n)}(k)=\varsigma_{n-1}\mu_{k}\left(v^{(n-1)}_{k}\right)^{p-1}, \qqd k \in E,\\
  v^{(n)}_{0}=0,
  \end{cases}
\end{equation*}
 where $\varsigma_{n-1}=\|w^{(n)}\|^{1-p}_{\mu,p}\|w^{(n-1)}\|^{p-1}_{\mu,p}$, and $$\lim_{n\rightarrow \infty}z_{n}=\lim_{n\rightarrow \infty}\varsigma_{n}=\varsigma.$$
 Here we omit the details of the proof.
\noindent\\[4mm]
\noindent\bf{\footnotesize Acknowledgements}\quad\rm
{\footnotesize This paper is based on a series of studies by Mu-Fa Chen. Heartfelt thanks are given to Professor Chen for his careful guidance and helpful suggestions.
 This work is supported in part by the National Nature Science Foundation of China (Grant No. 11771046), the project from the Ministry of Education in China.}\\[4mm]

\noindent{\bbb{References}}

\end{document}